\documentclass[11pt,reqno]{amsart}
\usepackage{amsmath, amssymb, amsthm}
\usepackage{url}
\usepackage[breaklinks]{hyperref}
\usepackage{array}
\renewcommand{\(}{\left\(}
\renewcommand{\)}{\right\)}
\renewcommand{\[}{\left\[}
\renewcommand{\]}{\right\]}

\numberwithin{equation}{section}
\theoremstyle{plain}
\newtheorem{theorem}{Theorem}[section]
\newtheorem{lemma}[theorem]{Lemma}

\newtheorem{definition}{Definition}[section]

\makeatletter
\def\proof{\@ifnextchar[{\@oproof}{\@nproof}}
\def\@oproof[#1][#2]{\trivlist\item[\hskip\labelsep\textit{#2 Proof of\
		#1.}~]\ignorespaces}
\def\@nproof{\trivlist\item[\hskip\labelsep\textit{Proof.}~]\ignorespaces}

\makeatother

\newenvironment{proof-alt}
{\vskip 0.15in \par\noindent{\it Proof of Proposition \ref{HSS}.}\hskip 0.5em\ignorespaces}
{\hfill $\Box$\par\medskip}

\begin{document}
	\title{Arithmetic properties of $5$-regular partitions into distinct parts}
	
	\author{Nayandeep Deka Baruah}
	\address{Nayandeep Deka Baruah, Department of Mathematical Sciences, Tezpur University, Napaam, Assam, 784028, India.}
	\email{nayan@tezu.ernet.in, nayandeeptezu@gmail.com}
	
	\author{Abhishek Sarma}
	\address{Abhishek Sarma, Department of Mathematical Sciences, Tezpur University, Napaam, Assam, 784028, India.}
	\email{abhitezu002@gmail.com}
	
	\thanks{2020 \textit{Mathematics Subject Classification.} 05A17, 11P81, 11P83, 11B37.\\
		\textit{Keywords and phrases. $\ell$-Regular partition, congruence, modular forms.}}
	
	\begin{abstract}
			A partition is said to be $\ell$-regular if none of its parts is a multiple of $\ell$.	Let $b^\prime_5(n)$ denote the number of 5-regular partitions into distinct parts (equivalently, into odd parts) of $n$. This function has also close connections to representation theory and combinatorics. In this paper, we study arithmetic properties of $b^\prime_5(n)$. We provide full characterization of the parity of $b^\prime_5(2n+1)$, present several congruences modulo 4, and prove that the generating function of the sequence $(b^\prime_5(5n+1))$ is lacunary modulo any arbitrary positive powers of 5.
	\end{abstract}
	\maketitle
	\section{Introduction}
	A partition $\lambda=(\lambda_1, \lambda_2, \ldots, \lambda_k)$ of a positive integer $n$ is a non-increasing sequence of positive integers, $\lambda_1\geq \lambda_2\geq \cdots \geq \lambda_k$ such that $\lambda_1+\lambda_2+\cdots+\lambda_k=n$. Here each $\lambda_{i}$ is called a part of the partition. For example, 2+2+1 is a partition of 5. The partition function $p(n)$ counts the number of partitions of $n$. 	
	
	A partition is said to be $\ell$-regular if none of its parts is a multiple of $\ell$.	For example, 3+3+2+1 is a 4-regular partition of 9. Let $b_{\ell}(n)$ denote the number of $\ell$-regular partitions of $n$. Then, with the convention, $b_{\ell}(0)=1$, the generating function of $b_{\ell}(n)$ is given by
	\begin{align*}
		\sum_{n=0}^{\infty} b_{\ell}(n)q^n=\frac{(q^{\ell};q^{\ell})_{\infty}}{(q;q)_{\infty}},
	\end{align*}
	where for complex numbers $a$ and $q$ with $|q|<1$, 
	\begin{align*}(a;q)_\infty:=\prod_{k=0}^\infty(1-aq^k).
	\end{align*}
	The $\ell$-regular partition function $b_\ell(n)$ has been studied quite extensively in the recent past by various authors. It is to be noted that the function $b_\ell(n)$ for prime $\ell$ gives the number of irreducible $\ell$-modular representation of the symmetric group $S_n$ \cite{james-kerber}. One can see the following non-exhaustive list of papers for various works related to $b_\ell(n)$ for $\ell\geq3$ (note that $b_2(n)$ counts the number of partitions of $n$ into odd parts, which is, by Euler's theorem, equal to the number of partitions of $n$ into distinct parts); arranged in alphabetical order of the first authors: 
	
	Ahmed and Baruah \cite{ahmed-baruah}, Alladi \cite{alladi}, Andrews, Hirschhorn, and Sellers \cite{andrews-hirschhorn-sellers}, Ballantine and Merca \cite{ballantine-merca}, Barman, Singh, and Singh \cite{Barman-AnnComb24}, Baruah and Das \cite{baruahdask},  Calkin, Drake, James, Law, Lee, Penniston, and Radder \cite{calkin-etal}, Carlson and Webb \cite{carlson-webb},  Cui and Gu \cite{cui-gu1,cui-gu2,cui-gu3, cui-gu4},  Dai \cite{dai}, Dai and Yan \cite{dai-yan}, Dandurand and Penniston \cite{dandurand-penniston}, Furcy and Penniston \cite{furcy-penniston}, Gordon and Ono \cite{gordon-ono}, Granville and Ono \cite{granville-ono}, Hirschhorn and Sellers \cite{hirschhorn-sellers},  Hou, Sun, and Zhang \cite{hou}, Iwata \cite{iwata},  Keith \cite{keith},  Keith and Zanello \cite{keith-zanello1, keith-zanello2}, Lin and Wang \cite{lin-wang},   Lovejoy \cite{lovejoy}, Lovejoy and Penniston \cite{lovejoy-penniston}, Mestrige \cite{mestrige}, Ono and Penniston \cite{ono-penniston1, ono-penniston2}, Penniston \cite{penniston1, penniston2, penniston3}, Singh and Barman \cite{singh-barman1, singh-barman2},   Singh, Singh, and Barman \cite{SinghEtAl-ActaArith}, Wang \cite{wang}, Webb \cite{webb}, Xia \cite{xia}, Xia and Yao \cite{xia-yao1, xia-yao2}, Yao \cite{yao1, yao2}, Zhao, Jin, and  Yao \cite{zhao-jin-yao}.  
	
	Let $b^\prime_\ell(n)$ count the number of $\ell$-regular partitions into distinct parts of $n$. For example, $b^\prime_{5}(10)=7$ and the relevant 7 partitions of 10 are $9+1$, $8+2$, $7+3$, $7+2+1$, $6+4$, $6+3+1$, and $4+3+2+1$. It is clear that $b^\prime_{\ell}(n)$ also counts the number of $\ell$-regular partitions with odd parts of $n$. With the convention that $b^\prime_{\ell}(0)=1$, the  generating function of $b^\prime_{\ell}(n)$ is given by
	\begin{align}\label{dk}
		\sum_{n=0}^\infty b^\prime_{\ell}(n)q^n&=\frac{(-q;q)_\infty}{(-q^\ell;q^\ell)_\infty}.
	\end{align}
	
	Note that, $b^\prime_2(n)$ counts the number of partitions of $n$ into  distinct odd parts, which, in fact, is equal to the number of self-conjugate partitions of $n$. The function  $b^\prime_2(n)$  has been well-studied. There are certain known results on $b^\prime_{\ell}(n)$ for $\ell\geq3$. For primes $\ell\geq 3$ and an integer $r$ with $1\leq r\leq p-1$ such that $24r+1$ is a quadratic nonresidue modulo $\ell$, Sellers \cite{sellers} proved that, for all nonnegative integers $n$,
	\begin{align}\label{sellers-cong}b^\prime_\ell(\ell n+r)\equiv 0~(\textup{mod}~ 2).
	\end{align}
	For a given prime $\ell\geq5$, Cui and Gu \cite[p. 523]{cui-gu1} showed that 
	\begin{align}\label{cui-gu-cong}b^\prime_\ell\left(\ell n+\dfrac{\ell^2-1}{24}\right)\equiv b_\ell(n)~(\textup{mod}~ 2).
	\end{align}
	Therefore, congruences modulo 2 of $b^\prime_\ell(n)$ may be studied from those of  $b_\ell(n)$. Thus,  many results on congruences modulo 2 for $b^\prime_\ell(n)$ can be derived from the results in papers on $b_\ell(n)$ that we mentioned earlier. Recently, Iwata \cite{iwata} found some congruences modulo 2 for  $b^\prime_\ell(n)$ for $\ell=9,25,41$, and $45$ by using modular forms.
	
	In this paper, we study the arithmetic properties of the function $b^\prime_5(n)$, which counts the number of 5-regular partitions into distinct parts of $n$. Setting $\ell=5$ in \eqref{dk}, we have
	\begin{align}\label{gf-bprime-5}
		\sum_{n=0}^\infty b^\prime_5(n)q^n&=\dfrac{(-q;q)_\infty}{(-q^5;q^5)_\infty}\notag\\
		&=1+q+q^2+2q^3+2q^4+2q^5+3q^6+4q^7+4q^8+6q^9+7q^{10}++8q^{11}\notag\\
		&\quad+10q^{12}+12q^{13}+14q^{14}+16q^{15}+19q^{16}+22q^{17}+26q^{18}+\cdots.
	\end{align}
	The sequence $(b^\prime_5(n))$ is A096938 in \cite{OEIS}. Other interpretations of this sequence are also discussed there. The function $b^\prime_5(n)$ is also related to representation theory and studied from that point of view by Andrews, Bessenrodt, and Olsson \cite{andrews-bessenrodt-olsson-archmath} (see also Andrews, Bessenrodt, and Olsson \cite{andrews-bessenrodt-olsson-tams}). Very recently, Ballantine and Feigon \cite{ballantine-feigon} gave a new combinatorial interpretation
	of $b^\prime_5(n)$. There are a few known arithmetic properties of $b^\prime_5(n)$ as well. It follows from \eqref{sellers-cong} that
	\begin{align*}b^\prime_5(5n+3)\equiv b^\prime_5(5n+4)\equiv0~(\textup{mod}~2).\end{align*}
	Many results on congruences modulo 2 for $b^\prime_5(n)$ can also be derived  from  \eqref{cui-gu-cong} and the corresponding work on $b_5(n)$. In particular, see \cite{calkin-etal},  \cite{cui-gu1}, and \cite{hirschhorn-sellers} for results on $b_5(n)$  modulo 2. 
	
	In this paper, we prove several new arithmetic results on $b^\prime_5(n)$. We state our results in the following theorems.
	
	The following theorem gives a complete characterization of the parity of $b^\prime_5(2n+1)$.
	\begin{theorem}\label{character2}For all $n \geq 0$, 
		\begin{align}\label{mod2}b^\prime_5(2n+1)\equiv \begin{cases}
				1\pmod2, & \text{if  $n=15k^2-5k$ for $k\in\mathbb{Z}$,}\\
				0\pmod2, & \text{Otherwise.} 
			\end{cases}
		\end{align}
	\end{theorem}
	Some congruences modulo 4 for $b^\prime_5(n)$ are given in the next theorem. To state the theorem, we require the Legendre symbol, which is defined for a prime $p\geq3$ by
	\begin{align*}
		\left(\dfrac{a}{p}\right)_L:=\begin{cases}\quad1,\quad \text{if $a$ is a quadratic residue modulo $p$ and $p\nmid a$,}\\\quad 0,\quad \text{if $p\mid a$,}\\~-1,\quad \text{if $a$ is a quadratic nonresidue modulo $p$.}
		\end{cases}
	\end{align*}
	
	\begin{theorem}\label{thm1}
		Let $p$ $(\geq5)$ be a prime such that $\left(\frac{3}{p}\right)_L\neq\left(\frac{-5}{p}\right)_L$. For all $n\geq0$ and $\alpha\geq0$, we have
		\begin{align}
			b^\prime_{5}(20n+j)&\equiv 0 \pmod4,~where ~ j \in  \{7, 15\},\label{cong-d5-20n}\\
			b^\prime_{5}(100n+j)&\equiv 0 \pmod4,~where ~j \in \{11, 31\},\label{cong-d5-100n}\\
			b^\prime_{5}\left(4\cdot p^{2\alpha}(5n+j)+\frac{17\cdot p^{2\alpha}+1}{6}\right)&\equiv 0 \pmod4,~where ~j \in \{1, 3\},\label{cong-d5-4-1-infty}\\
			b^\prime_{5}\left(4\cdot p^{2\alpha+1}(pn+j)+\frac{17\cdot p^{2\alpha+2}+1}{6}\right)&\equiv 0 \pmod4,~where ~j  \in\{1, 2, \ldots, p-1\}. \label{cong-d5-4-2-infty}
		\end{align}
	\end{theorem}
	
	In the next theorem, we state the exact generating functions of $b^\prime_5(5n+1)$ and $b^\prime_5(25n+21)$  in terms of $q$-products,  from which an internal congruence modulo 5 is also derived.	
	\begin{theorem}\label{exact}
		We have
		\begin{align}\label{exact1}
			\sum_{n=0}^{\infty} b^\prime_{5}(5n+1)q^n&=\frac{(q^2;q^2)_\infty (q^5;q^5)_\infty^3}{(q;q)_\infty^3 (q^{10};q^{10})_\infty}\\\intertext{and}
			\label{exact2}
			\sum_{n=0}^\infty b^\prime_5(25n+21)q^n&=\dfrac{(q;q)_\infty(q^{10};q^{10})_\infty^3}{(q^2;q^2)_\infty^3(q^5;q^5)_\infty}+40\dfrac{(q^2;q^2)_\infty^4
				(q^5;q^5)_\infty^4}{(q;q)_\infty^8}\notag\\
			&\quad+500q\dfrac{(q^2;q^2)_\infty^4
				(q^5;q^5)_\infty^{10}}{(q;q)_\infty^{14}}.
		\end{align}
		Furthermore, for all integers $\alpha\geq0$, we have
		\begin{align}\label{inffam}
			b^\prime_{5}(5n+1)\equiv b^\prime_{5}\left(5^{2\alpha+1}n+\frac{5^{2\alpha+1}+1}{6}\right)\pmod5.
		\end{align}
	\end{theorem}
	
	With the aid of \eqref{exact1}, we also study the distribution of $b^\prime_5(5n+1)$ modulo powers of 5.	
	
	Given any integral power series $F(q):=\sum_{n=0}^{\infty}a(n)q^n$ and $0\leq r<M$, define
	\begin{align*}
		\delta_r(F,M;X)	:= \frac{\#\{ n \leq X : a(n)\equiv r \pmod{M}\}}{X}.
	\end{align*}
	The series $F(q)$ is called lacunary modulo $M$ if
	\begin{align*}
		\lim_{X\to\infty}\delta_0(F,M;X) = 1.
	\end{align*} 
	
	We prove the following theorem which implies that $\sum_{n=0}^{\infty}b^\prime_5(5n+1)q^n$ is lacunary modulo arbitrary positive powers of 5.		
	\begin{theorem}\label{density}
		Let $k$ be a positive integer. Then
		\begin{align*}
			\lim_{X\to\infty} \frac{\#\{0\leq n \leq X : b^\prime_5 (5n+1)\equiv 0 \pmod{5^k}\}}{X} = 1.
		\end{align*}
	\end{theorem}
	
	In Sections \ref{sec3}--\ref{sec7}, we prove Theorems \ref{character2}--\ref{density}, respectively. We use $t$-dissections of certain $q$-products and the theory of modular forms in our proofs. The necessary background material and useful preliminary lemmas are given in the corresponding section.
	
	\section{Proof of Theorem \ref{character2}}\label{sec3}
	From \eqref{gf-bprime-5}, we have
	\begin{align}\label{b-mod2-1}
		\sum_{n=0}^\infty b^\prime_5(n)q^n=\dfrac{(-q;q)_\infty}{(-q^5;q^5)_\infty}=\dfrac{(q^2;q^2)_\infty(q^5;q^5)_\infty}{(q;q)_\infty(q^{10};q^{10})_\infty}.
	\end{align}
	
	From \cite[Theorem 2.1]{hirschhorn-sellers}, we recall that
	\begin{align}\label{b-mod2-2}
		\dfrac{(q^5;q^5)_\infty}{(q;q)_\infty}=\dfrac{(q^8;q^8)_\infty(q^{20};q^{20})_\infty^2}{(q^2;q^2)_\infty^2(q^{40};q^{40})_\infty}+q\dfrac{(q^4;q^4)_\infty^3(q^{10};q^{10})_\infty(q^{40};q^{40})_\infty}{(q^2;q^2)_\infty^3(q^8;q^8)_\infty(q^{20};q^{20})_\infty}.
	\end{align}	
	Employing \eqref{b-mod2-2} in \eqref{b-mod2-1}, we have
	\begin{align}\label{b-mod2-3}
		\sum_{n=0}^\infty b^\prime_5(n)q^n=\dfrac{(q^8;q^8)_\infty(q^{20};q^{20})_\infty^2}{(q^2;q^2)_\infty(q^{10};q^{10})_\infty(q^{40};q^{40})_\infty}+q\dfrac{(q^4;q^4)_\infty^3(q^{40};q^{40})_\infty}{(q^2;q^2)_\infty^2(q^8;q^8)_\infty(q^{20};q^{20})_\infty}.
	\end{align}
	Extracting the terms involving $q^{2n+1}$ from both sides of the above, dividing by $q$, and then replacing $q^2$ by $q$ in the resulting identity, we find that	
	\begin{align}\label{b-mod2-4}
		\sum_{n=0}^\infty b^\prime_5(2n+1)q^n=\dfrac{(q^2;q^2)_\infty^3(q^{20};q^{20})_\infty}{(q;q)_\infty^2(q^4;q^4)_\infty(q^{10};q^{10})_\infty}.
	\end{align}
	
	Now, it is easy to see from the binomial theorem that, for all positive integers $j$  
	\begin{align}\label{binom-mod2}
		(q^j;q^j)_\infty^2\equiv (q^{2j};q^{2j})_\infty\pmod2.
	\end{align}	
	Employing \eqref{binom-mod2} in \eqref{b-mod2-4}, we have
	\begin{align}\label{b-mod2-5}
		\sum_{n=0}^\infty b^\prime_5(2n+1)q^n\equiv(q^{10};q^{10})_\infty\pmod2.
	\end{align}
	
	However, Euler's famous pentagonal number theorem (see \cite[p. 11, Corollary 1.7]{andrews}) states that 
	\begin{align}\label{epnt}
		(q;q)_\infty=\sum_{k=-\infty}^\infty (-1)^kq^{k(3k-1)/2}.
	\end{align}
	Employing  the above, with $q$ replaced by $q^{10}$, in \eqref{b-mod2-5}, we see that 
	\begin{align*}
		\sum_{n=0}^\infty b^\prime_5(2n+1)q^n\equiv\sum_{k=-\infty}^\infty q^{15k^2-5k}\pmod2,
	\end{align*}
	from which we readily arrive at \eqref{mod2} to complete the proof.

	\section{Proof of Theorem \ref{thm1}}\label{sec4}
	
	To prove Theorem \ref{thm1}, we need some well-known results on $t$-dissections, where a $t$-dissection of a power series $A(q)$ in $q$ is given by
	$$A(q) =\sum_{j=0}^{t-1}q^jA_j(q^t),$$
	where $A_j$'s are power series in $q$.
	
	We also need results on Ramanujan's theta function $f(a,b)$ defined by
	\begin{align*}
		f(a,b) = \sum_{k=-\infty}^{\infty} a^{\frac{k(k+1)}{2}}b^{\frac{k(k-1)}{2}}, \quad |ab| < 1.
	\end{align*}
	Jacobi's triple product identity \cite[p. 35, Entry 19]{bcb3} is
	\begin{align*}
		f(a,b) := (-a;ab)_\infty(-b;ab)_\infty(ab;ab)_\infty,
	\end{align*}
	from which it follows that
	\begin{align*}
		f(-q,-q^2)=(q;q)_\infty,
	\end{align*}
	which is equivalent to \eqref{epnt}. 
	
	Now we state some useful lemmas. The first two lemmas below give 2- and 5-dissections of $(q;q)_\infty^2$ and $(q;q)_\infty$, respectively,  whereas the last two lemmas state $p$-dissections of $(q;q)_\infty$ and $(q;q)_\infty^3$ for primes $p\geq3$. 
	\begin{lemma}\textup{(\cite[From Eqs. (3.6.1) and (3.6.2)]{spirit})}
		\begin{align}
			(q;q)_\infty^2&= \frac{(q^2;q^2)_\infty (q^8;q^8)_\infty^5}{(q^4;q^4)_\infty^2(q^{16};q^{16})_\infty^2} - 2 q \frac{(q^2;q^2)_\infty (q^{16};q^{16})_\infty^2}{(q^8;q^8)_\infty}.\label{disf1^2}
		\end{align}
	\end{lemma}

	\begin{lemma}\label{q-5R}\textup{(\cite[p. 165]{spirit})} Let
		\begin{align*}
			R(q):= \dfrac{(q;q^5)_\infty(q^4;q^5)_\infty}{(q^2;q^5)_\infty(q^3;q^5)_\infty}.
		\end{align*}Then
		\begin{align}\label{disf_1}
			(q;q)_\infty&=(q^{25};q^{25})_\infty\left(\frac{1}{R(q^5)}-q-q^2R(q^5) \right).
		\end{align}
	\end{lemma}
	
	\begin{lemma}\textup{(\cite[Theorem 2.2]{cui-gu1})}\label{pdis}
		For a prime $p \geq 5$, we have 
		\begin{align}
			(q;q)_\infty&= (-1)^{\frac{\pm p-1}{6}} q^{\frac{p^2 - 1}{24}} (q^{p^2};q^{p^2})_\infty \notag\\
			&\quad+ \sum_{\begin{array}{c}
					k=-\frac{p-1}{2} \\
					k\neq\frac{\pm p-1}{6}
			\end{array}}^{\frac{p-1}{2}} (-1)^k q^{\frac{3k^2 + k}{2}} f\Bigg(-q^{\frac{3p^2+(6k+1)p}{2}},-q^{\frac{3p^2-(6k+1)p}{2}}\Bigg),\label{disf1}
		\end{align} 
		where 
		\begin{center}
			$\dfrac{\pm p-1}{6}=$ 
			$\begin{cases}
				\dfrac{p-1}{6}, & if ~p\equiv~ 1~ \pmod{6},\\
				\dfrac{-p-1}{6}, & if ~p\equiv~ -1~ \pmod{6}.
			\end{cases}$
		\end{center}
		Furthermore, for $\dfrac{-(p-1)}{2} \leq k \leq \dfrac{p-1}{2}$ and $k \neq \dfrac{\pm p-1}{6}$,
		\begin{align*}
			\dfrac{3k^2 + k}{2} \not\equiv \dfrac{p^2 - 1}{24} \pmod{p}.
		\end{align*}
	\end{lemma}
	
	\begin{lemma}\textup{(\cite[Lemma 2.3]{ahmed-baruah})}
		If $p\geq3$ is a prime, then
		\begin{align}\label{f13dis}
			(q;q)_\infty^3&=\sum_{\begin{array}{c}
					k=0 \\
					k\neq\frac{p-1}{2}
			\end{array}}^{p-1} (-1)^{k} q^{\frac{k(k+1)}{2}} \sum_{n=0}^{\infty} (-1)^n (2pn+2k+1) q^{pn\cdot\frac{pn+2k+1}{2}} \notag\\
			&\quad+ p (-1)^{\frac{p-1}{2}} q^{\frac{p^2-1}{8}} (q^{p^2};q^{p^2})_\infty^3.
		\end{align}
		Furthermore, for $0\leq k \leq p-1$ and $k\neq\frac{p-1}{2}$, 		\begin{align*}
			\frac{k^2+k}{2} \not\equiv \frac{p^2-1}{8} \pmod{p}.
		\end{align*}		
	\end{lemma}
	
	Now we are in a position to prove Theorem \ref{thm1}. 
	
	\noindent\emph{Proofs of \eqref{cong-d5-20n} and \eqref{cong-d5-100n}}.
	By the binomial theorem, for all positive integers $j$,  
	\begin{align}\label{binom-mod4}
		(q^j;q^j)_\infty^4\equiv (q^{2j};q^{2j})_\infty^2(\textup{mod}~4).
	\end{align}	
	Employing \eqref{binom-mod4} in \eqref{b-mod2-4}, and then using \eqref{disf1^2}, we have
	\begin{align*}
		\sum_{n\geq0}b^\prime_{5}(2n+1)q^n&\equiv \frac{(q;q)_\infty^2(q^2;q^2)_\infty(q^{20};q^{20})_\infty}{(q^4;q^4)_\infty (q^{10};q^{10})_\infty}\\
		&\equiv \frac{(q^2;q^2)_\infty^2(q^8;q^8)_\infty(q^{20};q^{20})_\infty}{(q^4;q^4)_\infty^3 (q^{10};q^{10})_\infty}\\
		&\quad-2q\frac{(q^2;q^2)_\infty^2(q^{16};q^{16})_\infty^2
			(q^{20};q^{20})_\infty}{(q^4;q^4)_\infty (q^8;q^8)_\infty(q^{10};q^{10})_\infty}\pmod4.
	\end{align*}
	Extracting the terms involving $q^{2n+1}$, we find that
	\begin{align*}
		\sum_{n\geq0}b^\prime_{5}(4n+3)q^n\equiv 2 \frac{(q;q)_\infty^2(q^8;q^8)_\infty^2(q^{10};q^{10})_\infty}{(q^2;q^2)_\infty(q^4;q^4)_\infty(q^5;q^5)_\infty}\pmod4.
	\end{align*}
	Applying \eqref{binom-mod2} in the above, we have
	\begin{align}
		\sum_{n\geq0}b^\prime_{5}(4n+3)q^n\equiv 2 (q^4;q^4)_\infty^3(q^5;q^5)_\infty\pmod4,\label{cong-d5-4n+3}
	\end{align}
	which, by \eqref{disf_1}, yields
	\begin{align}
		\sum_{n\geq0}b^\prime_{5}(4n+3)q^n&\equiv 2 (q^5;q^5)_\infty(q^{100};q^{100})_\infty^3
		\Bigg(\frac{1}{R(q^{20})^3}+\dfrac{q^4}{R(q^{20})^2}+q^{12}\notag\\
		&\quad+q^{20}R(q^{20})^2-\dfrac{q^{24}}{R(q^{20})^3}\Bigg)\pmod4.\label{cong-d5-4n+3-b}
	\end{align}	
	Equating the coefficients of $q^{5n+1}$ and $q^{5n+3}$ from both sides of the above, we arrive at
	\begin{align*}
		b^\prime_{5}(20n+7)\equiv b^\prime_{5}(20n+15)\equiv 0 \pmod4,
	\end{align*}
	which is  \eqref{cong-d5-20n}.
	
	Next, extracting the coefficients of $q^{5n+2}$ from both sides of \eqref{cong-d5-4n+3-b}, we find that
	\begin{align*}
		\sum_{n\geq0}b^\prime_{5}(20n+11)q^n\equiv2q^2(q;q)_\infty(q^{20};q^{20})_\infty^3\pmod4. 
	\end{align*} 
	Employing \eqref{disf_1} in the above and then equating the coefficients of $q^{5n}$ and $q^{5n+1}$, we obtain
	\begin{align*}
		b^\prime_{5}(100n+11)\equiv b^\prime_{5}(100n+31)\equiv0\pmod4,
	\end{align*}
	which is \eqref{cong-d5-100n}.
	
	\noindent \emph{Proofs of \eqref{cong-d5-4-1-infty} and \eqref{cong-d5-4-2-infty}}. 
	Employing \eqref{disf1} and \eqref{f13dis}, we rewrite \eqref{cong-d5-4n+3} as
	\begin{align*}
		&\sum_{n\geq0}b^\prime_{5}(4n+3)q^n\\
		&\equiv 2 \Bigg[\sum_{k\neq\frac{p-1}{2}, k=0}^{p-1} (-1)^{k} q^{4\cdot\frac{k(k+1)}{2}} \sum_{n=0}^{\infty} (-1)^n (2pn+2k+1) q^{4\cdot pn\cdot\frac{pn+2k+1}{2}} \\
		&\quad+ p (-1)^{\frac{p-1}{2}} q^{4\cdot \frac{p^2-1}{8}} (q^{4\cdot p^2};q^{4\cdot p^2})_\infty^3 \Bigg] \notag\\
		&\quad\times\Bigg[\sum_{k \neq \frac{\pm p-1}{6}, k = - \frac{p-1}{2} }^{\frac{p-1}{2}} (-1)^k q^{5\cdot\frac{3k^2 + k}{2}} f\Bigg(-q^{5\cdot\frac{3p^2+(6k+1)p}{2}},-q^{5\cdot\frac{3p^2-(6k+1)p}{2}}\Bigg)\notag\\
		&\quad+(-1)^{\frac{\pm p-1}{6}} q^{5\cdot\frac{p^2 - 1}{24}} (q^{5\cdot p^2};q^{5\cdot p^2})_\infty \Bigg]\pmod4.
	\end{align*}
	Consider the congruence
	\begin{align*}
		4\cdot\frac{k(k+1)}{2}+5\cdot\frac{3m^2 + m}{2}\equiv17\cdot\frac{p^2-1}{24}\pmod{p},
	\end{align*}
	where $0\leq k\leq p-1$ and $-\frac{(p-1)}{2}\leq m\leq \frac{p-1}{2}$. Since  the above congruence is equivalent to 
	\begin{align*}
		3(4k+2)^2+5(6m+1)^2\equiv0\pmod{p}
	\end{align*}
	and $\left(\frac{3}{p}\right)_L\neq\left(\frac{-5}{p}\right)_L$, it  turns out that the only possibilities of satisfying the above congruence are $k=\frac{p-1}{2}$ and $m=\frac{p-1}{6}$. So, extracting the terms involving $q^{pn+17\cdot\frac{p^2-1}{24}}$,  dividing both sides by $q^{17\cdot\frac{p^2-1}{24}}$, and then replacing $q^p$ by  $q$, we arrive at
	\begin{align*}
		\sum_{n\geq0}b^\prime_{5}\left(4pn+\frac{17\cdot p^2+1}{6}\right)q^n&\equiv 2p(q^{4p};q^{4p})_\infty^3(q^{5p};q^{5p})_\infty \pmod4.
	\end{align*}
	Extracting the terms involving $q^{pn}$, we have
	\begin{align*}
		\sum_{n\geq0}b^\prime_{5}\left(4p^2n+\frac{17\cdot p^2+1}{6}\right)q^n&\equiv 2p(q^{4};q^{4})_\infty^3(q^{5};q^{5})_\infty \pmod4.
	\end{align*}
	We now apply \eqref{disf1} and \eqref{f13dis} in the above and repeat the process $\alpha$ times to arrive at 
	\begin{align}
		\sum_{n\geq0}b^\prime_{5}\left(4p^{2\alpha}n+\frac{17\cdot p^{2\alpha}+1}{6}\right)q^n&\equiv 2p^{\alpha}f_{4}^3f_{5}\pmod4.\label{cong-d5-inftystep}
	\end{align}
	Using \eqref{f13dis} in the last step and extracting the terms involving  $q^{5n+1}$ and $q^{5n+3}$, we obtain \eqref{cong-d5-4-1-infty}.
	
	Again, extracting the terms involving $q^{pn+17\cdot\frac{p^2-1}{24}}$ from \eqref{cong-d5-inftystep},  dividing both sides by $q^{17\cdot\frac{p^2-1}{24}}$, and then replacing $q^p$ by  $q$, we arrive at
	\begin{align*}
		\sum_{n\geq0}b^\prime_{5}\left(4p^{2\alpha+1}n+\frac{17\cdot p^{2\alpha+2}+1}{6}\right)q^n&\equiv 2p^{\alpha+1}f_{4p}^3f_{5p}\pmod4.
	\end{align*}
	Comparing the coefficients of $q^{pn+r}$, where $r\in\{1, 2, \ldots, p-1\}$, we obtain \eqref{cong-d5-4-2-infty}. This completes the proof of Theorem \ref{thm1}.

\section{Proof of Theorem \ref{exact}}\label{sec6}
First we state some lemmas.
\begin{lemma} \textup{(\cite[Eqs. (7.4.9) and (7.4.14)]{spirit})} If $R(q)$ is as defined in Lemma \ref{q-5R}, then
\begin{align}\label{nm}
	\dfrac{1}{R^5(q)}-q^2R^5(q)&=\frac{(q;q)_\infty^6}{(q^5;q^5)_\infty^6}+11q\\\intertext{and}\label{eq:bcb}
	\frac{1}{(q;q)_\infty}&=\frac{(q^{25};q^{25})_\infty^5}{(q^5;q^5)_\infty^6}\Big(\frac{1}{R(q^5)^4}+\frac{q}{R(q^5)^3}+\frac{2q^2}{R(q^5)^2}+\frac{3q^3}{R(q^5)}+5q^4 \notag\\
	&\quad-3q^5R(q^5)+2q^6R(q^5)^2-q^7R(q^5)^3+q^8R(q^5)^4\Big).
\end{align}
\end{lemma}
\begin{lemma}\textup{(\cite[Eqs (2.6) and (2.7)]{baruahbegum})} We have
\begin{align}
	\frac{(q^5;q^5)_\infty}{(q^2;q^2)_\infty^2(q^{10};q^{10})_\infty}	 &= \frac{(q^5;q^5)_\infty^5}{(q;q)_\infty^4(q^{10};q^{10})_\infty^3}-4q \frac{(q^{10};q^{10})_\infty^2}{(q;q)_\infty^3(q^2;q^2)_\infty},\label{nm1}\\
	\frac{(q^2;q^2)_\infty^3(q^5;q^5)_\infty^2}{(q;q)_\infty^2(q^{10};q^{10})_\infty^2} &= \frac{(q^5;q^5)_\infty^5}{(q;q)_\infty(q^{10};q^{10})_\infty^3}+q \frac{(q^{10};q^{10})_\infty^2}{(q^2;q^2)_\infty}.\label{nm2}
\end{align}
\end{lemma}
\begin{lemma}\textup{(\cite[Lemma 1.3]{baruahbegum})} If $R(q)$ is as defined in Lemma \ref{q-5R}, then
\begin{align}\label{xy2}
	\dfrac{1}{R(q)R^2(q^2)}-q^2R(q)R^2(q^2)&=\frac{(q^2;q^2)_\infty(q^5;q^5)_\infty^5}{(q;q)_\infty(q^{10};q^{10})_\infty^5},\\		
	\label{x2y}
	\dfrac{R(q^2)}{R^2(q)}-\dfrac{R^2(q)}{R(q^2})&=4q\frac{(q;q)_\infty(q^{10};q^{10})_\infty^5}{(q^2;q^2)_\infty(q^5;q^5)_\infty^5},\\
	\label{x3y}\dfrac{1}{R^3(q)R(q^2)}+q^2R^3(q)R(q^2)&=\frac{(q^2;q^2)_\infty(q^5;q^5)_\infty^5}{(q;q)_\infty(q^{10};q^{10})_\infty^5} +4q^2\frac{(q;q)_\infty(q^{10};q^{10})_\infty^5}{(q^2;q^2)_\infty(q^5;q^5)_\infty^5}+2q,\\
	\label{xy3}\dfrac{R(q)}{R^3(q^2)}+q^2\dfrac{R^3(q^2)}{R(q)}&=\frac{(q^2;q^2)_\infty(q^5;q^5)_\infty^5}{(q;q)_\infty(q^{10};q^{10})_\infty^5} +4q^2\frac{(q;q)_\infty(q^{10};q^{10})_\infty^5}{(q^2;q^2)_\infty(q^5;q^5)_\infty^5}-2q.
\end{align}			
\end{lemma}

Now we prove Theorem \ref{exact} by establishing \eqref{exact1}--\eqref{inffam}.	

\noindent \emph{Proof of \eqref{exact1}}.
Employing \eqref{disf_1}, with $q$ replaced by $q^2$ and \eqref{eq:bcb} in \eqref{b-mod2-1}, and then extracting the terms involving $q^{5n+1}$, we find that
\begin{align*}
\sum_{n=0}^\infty b^\prime_5(5n+1)q^n&=\dfrac{(q^5;q^5)_\infty^5(q^{10};q^{10})_\infty}{(q;q)_\infty^5(q^2;q^2)_\infty}\Bigg(\dfrac{1}{R^3(q)R(q^2)}+q^2R^3(q)R(q^2)\\
&\quad-2q\left(\dfrac{R(q^2)}{R^2(q)}-\dfrac{R^2(q)}{R(q^2)}\right)-5q\Bigg),
\end{align*}
which, by \eqref{x2y} and \eqref{x3y}, yields
\begin{align*}
\sum_{n=0}^\infty b^\prime_5(5n+1)q^n&=\dfrac{(q^5;q^5)_\infty^{10}}{(q;q)_\infty^6(q^{10};q^{10})_\infty^4}
-3q\dfrac{(q^5;q^5)_\infty^5(q^{10};q^{10})_\infty}{(q;q)_\infty^5(q^2;q^2)_\infty}-4q^2\dfrac{(q^{10};q^{10})_\infty^6}{(q;q)_\infty^4(q^2;q^2)_\infty^2}\\
&=\left(\dfrac{(q^5;q^5)_\infty^{10}}{(q;q)_\infty^6(q^{10};q^{10})_\infty^4}
-4q\dfrac{(q^5;q^5)_\infty^5(q^{10};q^{10})_\infty}{(q;q)_\infty^5(q^2;q^2)_\infty}\right)\\
&\quad+q\left(\dfrac{(q^5;q^5)_\infty^5(q^{10};q^{10})_\infty}{(q;q)_\infty^5(q^2;q^2)_\infty}-4q\dfrac{(q^{10};q^{10})_\infty^6}{(q;q)_\infty^4(q^2;q^2)_\infty^2}\right).
\end{align*}
Employing \eqref{nm1} and \eqref{nm2} in the above, we have
\begin{align}\label{exact-a}
\sum_{n=0}^\infty b^\prime_5(5n+1)q^n&=\dfrac{(q^5;q^5)_\infty^6}{(q;q)_\infty^2(q^2;q^2)_\infty^2(q^{10};q^{10})_\infty^2}
+q\dfrac{(q^5;q^5)_\infty(q^{10};q^{10})_\infty^3}{(q;q)_\infty(q^2;q^2)_\infty^3}\notag\\
&=\dfrac{(q^2;q^2)_\infty(q^5;q^5)_\infty^3}{(q;q)_\infty^3(q^{10};q^{10})_\infty},
\end{align}
which proves \eqref{exact1}.

\noindent \emph{Proof of \eqref{exact2}}. With the aid of \eqref{nm1}, we may rewrite \eqref{exact-a} as
\begin{align*}
\sum_{n=0}^\infty b^\prime_5(5n+1)q^n&=\dfrac{(q;q)_\infty(q^{10};q^{10})_\infty}{(q^2;q^2)_\infty(q^5;q^5)_\infty}+4q\dfrac{(q^{10};q^{10})_\infty^4}{(q;q)_\infty^2(q^5;q^5)_\infty^2}.
\end{align*}
Employing \eqref{disf_1} and \eqref{eq:bcb} in the above identity, and then extracting the terms involving $q^{5n+4}$, we find that	
\begin{align*}
\sum_{n=0}^\infty b^\prime_5(25n+21)q^n&=\dfrac{(q^5;q^5)_\infty(q^{10};q^{10})_\infty^5}{(q;q)_\infty(q^2;q^2)_\infty^5}\Bigg(2\left(\dfrac{1}{R(q)R^2(q^2)}-q^2R(q)R^2(q^2)\right)\notag\\
&\quad-\left(\dfrac{R(q)}{R^3(q^2)}+q^2\dfrac{R^3(q^2)}{R(q)}\right)-5q\Bigg)\notag\\
&+20\dfrac{(q^2;q^2)_\infty^4(q^5;q^5)_\infty^{10}}{(q;q)_\infty^{14}}\left(2\left(\dfrac{1}{R^5(q)}-q^2R^5(q)\right)+3q\right).
\end{align*}
Using \eqref{nm}, \eqref{xy2},  and \eqref{xy3} in the above identity, we have
\begin{align*}
\sum_{n=0}^\infty b^\prime_5(25n+21)q^n&=\dfrac{(q^5;q^5)_\infty^6}{(q;q)_\infty^2(q^2;q^2)_\infty^4}-4q\dfrac{(q^5;q^5)_\infty
	(q^{10};q^{10})_\infty^5}{(q;q)_\infty(q^2;q^2)_\infty^5}\\
&\quad+q\left(\dfrac{(q^5;q^5)_\infty
	(q^{10};q^{10})_\infty^5}{(q;q)_\infty(q^2;q^2)_\infty^5}-4q\dfrac{(q^{10};q^{10})_\infty^{10}}{(q^2;q^2)_\infty^6(q^5;q^5)_\infty^4}\right)\\
&\quad+40\dfrac{(q^2;q^2)_\infty^4
	(q^5;q^5)_\infty^4}{(q;q)_\infty^8}+500q\dfrac{(q^2;q^2)_\infty^4
	(q^5;q^5)_\infty^{10}}{(q;q)_\infty^{14}}.
\end{align*}
With the aid of \eqref{nm1} and \eqref{nm2},  the above identity reduces to
\begin{align*}
\sum_{n=0}^\infty b^\prime_5(25n+21)q^n&=\dfrac{(q;q)_\infty(q^{10};q^{10})_\infty^3}{(q^2;q^2)_\infty^3(q^5;q^5)_\infty}+40\dfrac{(q^2;q^2)_\infty^4
	(q^5;q^5)_\infty^4}{(q;q)_\infty^8}\\
&\quad+500q\dfrac{(q^2;q^2)_\infty^4
	(q^5;q^5)_\infty^{10}}{(q;q)_\infty^{14}},
\end{align*}
which is \eqref{exact2}.

\noindent \emph{Proof of \eqref{inffam}}. By the binomial theorem, for all positive integers $j$, we have
\begin{align}\label{binom-mod5}
(q^j;q^j)_\infty^5\equiv (q^{5j};q^{5j})_\infty(\textup{mod}~5).
\end{align}	
Employing \eqref{binom-mod5} in \eqref{exact2}, we find that
\begin{align*}
\sum_{n=0}^\infty b^\prime_5(25n+21)q^n&\equiv\dfrac{(q;q)_\infty(q^{10};q^{10})_\infty^3}{(q^2;q^2)_\infty^3(q^5;q^5)_\infty}\\		&\equiv\dfrac{(q;q)_\infty(q^2;q^2)_\infty^2(q^{10};q^{10})_\infty^2}{(q^5;q^5)_\infty}~(\textup{mod}~5),
\end{align*}
which, with the help of \eqref{disf_1}, may be recast as
\begin{align*}
&\sum_{n=0}^\infty b^\prime_5(25n+21)q^n\\
&\equiv	\dfrac{(q^{10};q^{10})_\infty^2(q^{25};q^{25})_\infty(q^{50};q^{50})_\infty^2}{(q^5;q^5)_\infty}\left(\frac{1}{R(q^5)}-q-q^2R(q^5) \right)\\
&\quad\times\left(\frac{1}{R(q^{10})}-q^2-q^4R(q^{10}) \right)^2~(\textup{mod}~5).
\end{align*}
Extracting the terms involving $q^{5n}$ from both sides of the above, and then replacing $q^5$ by $q$, we obtain
\begin{align*}
&\sum_{n=0}^\infty b^\prime_5(125n+21)q^n\\
&\equiv	\dfrac{(q^2;q^2)_\infty^2(q^5;q^5)_\infty(q^{10};q^{10})_\infty^2}{(q;q)_\infty}\left(\frac{1}{R(q)R^2(q^2)}+q-q^2R(q)R^2(q^2) \right)~(\textup{mod}~5).
\end{align*} 
Employing \eqref{xy2} and \eqref{nm2} in the above, and then invoking \eqref{binom-mod5}, we find that
\begin{align*}
&\sum_{n=0}^\infty b^\prime_5(125n+21)q^n\equiv	\dfrac{(q^2;q^2)_\infty^6(q^5;q^5)_\infty^3}{(q;q)_\infty^3(q^{10};q^{10})_\infty^2}\equiv	\dfrac{(q^2;q^2)_\infty(q^5;q^5)_\infty^3}{(q;q)_\infty^3(q^{10};q^{10})_\infty}~(\textup{mod}~5).
\end{align*}
From the above congruence and \eqref{exact1}, it follows that	
\begin{align*}
b^\prime_5(5n+1)&\equiv b^\prime_{5}(125n+21)~(\textup{mod}~5),
\end{align*} which by iteration yields \eqref{inffam}. 

\section{Proof of Theorem \ref{density}}\label{sec7}

Before proving Theorem \ref{density}, we recall some useful background material on modular forms. Let $\mathbb{H}$ denote the complex upper half-plane. We define the following matrix groups:
\begin{align*}
\textup{SL}_2(\mathbb{Z}) &:=\left\{\begin{bmatrix}
	a  &  b \\
	c  &  d      
\end{bmatrix}: a, b, c, d \in \mathbb{Z}, ad-bc=1
\right\},\\
\Gamma_{0}(N) &:=\left\{
\begin{bmatrix}
	a  &  b \\
	c  &  d      
\end{bmatrix} \in \Gamma : c\equiv~0\pmod N \right\},\\
\Gamma_{1}(N) &:=\left\{
\begin{bmatrix}
	a  &  b \\
	c  &  d      
\end{bmatrix} \in \Gamma_0(N) : a\equiv~d\equiv~1\pmod N \right\},\\\intertext{and} 
\Gamma(N) &:=\left\{
\begin{bmatrix}
	a  &  b \\
	c  &  d      
\end{bmatrix} \in \textup{SL}_2(\mathbb{Z}) : a\equiv d\equiv 1\pmod N,~ \textup{and}~b\equiv c\equiv 0\pmod N\right\},
\end{align*}
where  $N$ is a positive integer.

A subgroup $\Gamma$ of $\textup{SL}_2(\mathbb{Z})$ satisfying $\Gamma(N)\subseteq\Gamma$ for some $N$ is called a congruence subgroup and the smallest such $N$ is called the level of $\Gamma$. The group
$$\textup{GL}^{+}_{2}(\mathbb{R}) :=\left\{
\begin{bmatrix}
a  &  b \\
c  &  d      
\end{bmatrix} : a, b, c, d \in \mathbb{R}~\textup{and}~ad-bc>0\right\}$$
acts on $\mathbb{H}$ by $\begin{bmatrix}
a  &  b \\
c  &  d      
\end{bmatrix}z = \dfrac{az+b}{cz+d}$. We will identify $\infty$ with $\dfrac{1}{0}$. We also define $\begin{bmatrix}
a  &  b \\
c  &  d
\end{bmatrix}\dfrac{r}{s} = \dfrac{ar+bs}{cr+ds}$, where $\dfrac{r}{s}\in\mathbb{Q}\cup\{\infty\}$. This will give an action of $\textup{GL}^{+}_{2}(\mathbb{R})$ on the extended upper half-plane $\mathbb{H^{*}}=\mathbb{H}\cup\mathbb{Q}\cup\{\infty\}$. If $\Gamma$ is a congruence subgroup of $\textup{SL}_2(\mathbb{Z})$, then a cusp of $\Gamma$ is an equivalence class in $\mathbb{Q}\cup\{\infty\}$.

For a positive integer $\ell$, let $M_\ell(\Gamma_1(N)))$ denote the complex vector space of modular forms of weight $\ell$ with respect to $\Gamma_1(N)$. 
\begin{definition}\cite[Definition 1.15]{OnoModularity}
If $\chi$ is a Dirichlet character modulo $N$, then a form $f(x)\in M_\ell(\Gamma_1(N))$ has Nebentypus character $\chi$ if 
$$f\left( \frac{az+b}{cz+d}\right)=\chi(d)(cz+d)^{\ell}f(z)$$ for all $z\in \mathbb{H}$ and all $\begin{bmatrix}
	a  &  b \\
	c  &  d      
\end{bmatrix}\in \Gamma_0(N)$. The space of such modular forms is denoted by $M_\ell(\Gamma_0(N),\chi)$.
\end{definition}

Now recall that the Dedekind's eta-function $\eta(z)$ is defined by
\begin{align*}
\eta(z):=q^{1/24}(q;q)_{\infty},
\end{align*}
where $q:=e^{2\pi iz}$ and $z\in \mathbb{H}$. A function $f(z)$ is called an eta-quotient if it is of the form
\begin{align*}
f(z)=\prod_{\delta|N}\eta(\delta z)^{r_\delta},
\end{align*}
where $N$ is a positive integer and $r_{\delta}$ is an integer. 

Next, we recall three theorems from \cite[p. 18]{OnoModularity} which will be used to prove Theorem \ref{density}.	
\begin{theorem}\cite[Theorem 1.64]{OnoModularity}\label{thm:ono1} If $f(z)=\prod_{\delta|N}\eta(\delta z)^{r_\delta}$ 
is an eta-quotient with $\ell=\dfrac{1}{2}\sum_{\delta|N}r_{\delta}\in \mathbb{Z}$, 
with
$$\sum_{\delta|N} \delta r_{\delta}\equiv 0 \pmod{24}$$ and
$$\sum_{\delta|N} \frac{N}{\delta}r_{\delta}\equiv 0 \pmod{24},$$
then $f(z)$ satisfies $$f\left( \frac{az+b}{cz+d}\right)=\chi(d)(cz+d)^{\ell}f(z)$$
for every  $\begin{bmatrix}
	a  &  b \\
	c  &  d      
\end{bmatrix} \in \Gamma_0(N)$, where $$\chi(d):=\left(\frac{(-1)^{\ell} \prod_{\delta |N}\delta^{r_{\delta}}}{d}\right).$$ 
\end{theorem}
\begin{theorem}\cite[Theorem 1.65]{OnoModularity}\label{thm:ono2}
If $c$, $d$, and $N$ are positive integers such that $d|N$ and $\gcd(c, d)=1$, then the order of vanishing of $f(z)$ at the cusp $\dfrac{c}{d}$ 
is
\begin{align*}
	\dfrac{N}{24}\sum_{\delta|N}\dfrac{\gcd(d,\delta)^2r_{\delta}}{\gcd(d,\frac{N}{d})d\delta}.
\end{align*} 
\end{theorem}
Suppose that $f(z)$ is an eta-quotient satisfying the conditions of the above theorem. 
If $f(z)$ is holomorphic at all of the cusps of $\Gamma_0(N)$, 
then $f(z)\in M_{\ell}(\Gamma_0(N), \chi)$. The following result is due to Serre, which we state from \cite[p. 43]{OnoModularity}).
\begin{theorem}\label{thm:ono3}
If $f(z)\in M_\ell(\Gamma_0(N), \chi)$ has Fourier expansion
\begin{align*}
	f(z) = \sum_{n=0}^\infty c(n)q^n\in \mathbb{Z}[[q]],
\end{align*}
then for each positive integer $m$ there exists a constant $\alpha>0$ such that
\begin{align*}
	|\{n\leq X~:~c(n)\not\equiv0\pmod m\}|=\mathcal{O}\left(\frac{X}{(\log X)^\alpha}\right).
\end{align*}
\end{theorem}

Now we are in a position to prove Theorem \ref{density}.

Let
\begin{align*}
A(z) := \prod_{n=1}^{\infty} \frac{(1-q^{12n})^5}{(1-q^{60n})} = \frac{\eta^5(12z)}{\eta(60z)}.
\end{align*}
Then
\begin{align*}
A^{5^k}(z) = \frac{\eta^{5^{k+1}}(12z)}{\eta^{5^k}(60z)}.
\end{align*}

Set
\begin{align*}
B_k(z) := \frac{\eta(12z)\eta^3(30z)}{\eta^3(6z)\eta(60z)} A^{5^k}(z)=\frac{\eta^{5^{k+1}+1}(12z)\eta^3(30z)}{\eta^3(6z)\eta^{5^{k}+1}(60z)}.
\end{align*}

Working modulo $5^{k+1}$, we have
\begin{align}\label{bk}
B_k(z) \equiv \frac{\eta(12z)\eta^3(30z)}{\eta^3(6z)\eta(60z)} = q \frac{(q^{12};q^{12})_\infty (q^{30};q^{30})_\infty^3}{(q^6;q^6)_\infty^3(q^{60};q^{60})_\infty}.
\end{align}
From \eqref{exact1} and \eqref{bk}, we see that
\begin{align}
B_k(z) \equiv \sum_{n=0}^{\infty} b^\prime_{5}(5n+1)q^{6n+1} \pmod{5^{k+1}}.\label{6n+1}
\end{align}

Clearly, $B_k(z)$ is an eta-quotient with $N=360$. We now prove that $B_k(z)$ is a modular form for any positive integer $k$. We know that the cusps of $\Gamma_0(360)$ are given by fractions $\frac{c}{d}$, where $d|360$ and $\gcd(c,d)=1$. By Theorem \ref{thm:ono2}, we find that $B_k(z)$ is holomorphic at a cusp $\frac{c}{d}$ if and only if
\begin{align*}
L:= (5^{k+2}+5)\frac{\gcd(d,12)^2}{\gcd(d,60)^2}+6\frac{\gcd(d,30)^2}{\gcd(d,60)^2}-30\frac{\gcd(d,6)^2}{\gcd(d,60)^2}-5^k-1\geq0.
\end{align*}
We verify that the above inequality holds for all the divisors of 360. We illustrate this with the help of the following table.
\begin{center}
\begin{tabular}{ | m{9em} | m{9cm}| } 
	\hline
	$d$ & $(5^{k+2}+5)\frac{\gcd(d,12)^2}{\gcd(d,60)^2}+6\frac{\gcd(d,30)^2}{\gcd(d,60)^2}-30\frac{\gcd(d,6)^2}{\gcd(d,60)^2}-5^k-1$\\ 
	\hline
	1,2,3,6,9,18 & $24\cdot5^k-20$  \\ 
	\hline
	4,8,12,24,36,72 & 24$\cdot5^{k}-2$\\ 
	\hline
	5,10,15,30,45,90 & $\frac{1}{25} \left(5^{k+2}+5\right)-5^k+\frac{19}{5}$\\
	\hline
	20,40,60,120,180,360 & $\frac{1}{25} \left(5^{k+2}+5\right)-5^k+\frac{1}{5}$\\
	\hline
\end{tabular}
\end{center}
Using Theorem \ref{thm:ono1}, we find that the weight of $B_k(z)$ is $2\cdot5^k$. Further, the associated character for $B_k(z)$ is $\chi_1(\bullet)=\left(\frac{12^{4\cdot5^k}\cdot5^{2-5^k}}{\bullet}\right)$. Thus, $B_k(z) \in M_{2\cdot5^k}(\Gamma_0(360), \chi_1)$. Also, the Fourier coefficients of $B_k(z)$ are all integers. Hence, by Theorem \ref{thm:ono3}, the Fourier coefficients of $B_k(z)$ are almost always divisible by $5^k$. By \eqref{6n+1}, the same holds for $b^\prime_5(5n+1)$ and hence Theorem \ref{density} follows.
\subsection*{Acknowledgements} A. Sarma was partially supported by an institutional fellowship for doctoral research from Tezpur University, Assam, India. The author thanks the funding institution. The authors thank the anonymous referee for his/her careful reading of the manuscript and helpful comments.

\end{document}